\documentclass{amsart}
\title{Cluster algebras and Grassmannians of type $G_2$}
\author{Sachin Gautam}
\address{Department of Mathematics, Northeastern University, Boston MA 02115 (USA)}
\email{gautam.s@neu.edu}

\usepackage{amsmath, amssymb, hyperref}
\usepackage[all]{xy}

\newtheorem{thm}{Theorem}[section]
\newtheorem{prop}{Proposition}[section]
\newtheorem{rem}{Remark}[section]

\newtheorem{lem}{Lemma}[section]

\newtheorem{cor}{Corollary}[section]

\newcommand{\Lg}{\mathfrak{g}}

\newcommand{\C}{\mathbb{C}}
\newcommand{\Z}{\mathbb{Z}}
\newcommand{\A}{\mathcal{A}}
\newcommand{\N}{\mathbb{N}}
\newcommand{\F}{\mathcal{F}}
\newcommand{\Q}{\mathbb{Q}}
\usepackage{multirow}

\begin{document}
\maketitle

\begin{abstract}
We prove a conjecture of Geiss, Leclerc and Schr\"{o}er,  producing cluster algebra structures on multi-homogeneous coordinate ring of partial flag varieties, for the case $G_2$. As a consequence we sharpen the known fact that coordinate ring of the double Bruhat cell $G^{e,w_0}$ is an upper cluster algebra, by proving that it is a cluster algebra.
\end{abstract}

\tableofcontents

\section{Introduction}
The theory of cluster algebras was founded and developed by Sergey Fomin and Andrei Zelevinsky (see \cite{CA1}, \cite{CA3}, \cite{CA4}) in order to develop a framework for dual canonical basis. Since their first appearance, cluster algebras have found their place in several interesting branches of mathematics (e.g, Poisson geometry, integrable systems, Teichm\"{u}ller spaces, representations of finite dimensional algebras, tropical geometry,...).\\

The results of this note arose from an attempt to understand the cluster algebra structures coming from representations of the group of type $G_2$. The link between the representations of simple Lie group $G$ and cluster algebras has already been made in \cite{CA3} where authors give a combinatorial procedure to get a cluster algebra structure on the coordinate rings of double Bruhat cells, in particular on $G^{e,w_0}$ which is the one we will be focusing on in this note. (I also learnt from Prof. Andrei Zelevinsky that the computations of dual canonical basis for the case $G_2$ was a major motivation for defining cluster algebras).\\

The starting point of this project was the algebra of $N^-$ invariant functions on $G$ which serves as a model for all finite dimensional representations of $G$, and is equipped with cluster algebra structure. In the case of $G_2$ there are two fundamental representations and the algebra at hand is direct sum of finite dimensional representations parametrized by non-negative elements of two dimensional lattice (weight lattice, generated by fundamental weights $\omega_1, \omega_2$). 
\[
\A = \bigoplus_{m,n\in\Z_{\geq 0}} V_{m\omega_1 + n\omega_2}
\]

The following questions arise naturally in pursuit of this study:

\begin{itemize}
\item[(A)] For each $i = 1,2$ we have subalgebras $\A_i = \oplus_{n\in \Z_{\geq 0}} V_{n\omega_i}$ of $\A$. Is there a natural cluster algebra structure on $\A_i$?
\item[(B)] Is it possible to extract $\A_i$ from $\A$ by elementary operations (mutations, freezing of vertices, etc.) which will make the inclusion $\A_i \subset \A$ explicit?
\end{itemize}

In \cite{gls} authors develop a connection between cluster algebra structures on multi-homogeneous rings of partial flag varieties and representations of preprojective algebra of underlying Dynkin quiver, for simply laced cases. They also provide a purely combinatorial algorithm to extract these structures for non-simply laced cases (and hence giving a conjectural answer to (A)), but since the techniques of \cite{gls} are only applicable in simply-laced cases, the construction remained conjectural. The construction of  cluster algebra structures on $N$ (maximal unipotent group), using representation theory of preprojective algebras, was further extended to non simply-laced cases in \cite{demonet}. In loc. cit. author uses the folding process (of V. Kac and G. Lusztig) to reduce the problem to simply laced cases.\\

The aim of this note is to answer these two questions for case $G_2$. We produce cluster algebra structures on $\A_i$ starting from that of $\A$ by ``elementary operations". It turns out that the two structures are ``transpose" of each other (thus $\A_2$ is of type $G_2^{(1)}$ and $\A_1$ is of type $D_4^{(3)}$). Both these structure fall into the realm of \cite{giovanni}. In fact the statements of section \ref{propositions} are easy consequences of results of loc. cit. and are given here just for the sake of completeness. Next we prove that our construction agrees with that of \cite{gls}. This suggests an alternative method of constructing cluster algebra structures on multi-homogeneous coordinate rings of partial flag varieties, but I do not have the answer to this question in general. \\

The structure of this note is as follows. In section \ref{mainresults} we define the initial seed for each $\A_i$ and state the main theorem. We also observe that the statement of our main theorem (Theorem \ref{main})  implies conjectures 10.4 and 9.6 of \cite{gls} in $G_2$ case; and proves that upper cluster algebra structure on $G^{e,w_0}$ given in \cite{CA3} is in fact a cluster algebra structure.\\

Section \ref{review} is devoted to recall basic definitions from the theory of cluster algebras. In particular we recall the (upper) cluster algebra structure on $G^{e,w_0}$ as given in \cite{CA3}.\\

In section \ref{explanation} we explain how to get the initial seeds for $\A_i$ starting from initial seed of $\A$ and using some elementary operations of mutations and freezing of vertices. This allows us to use results of \cite{CA3} in our situation.\\

Sections \ref{proof}, \ref{propositions} and \ref{lemmas} are devoted to the proof of Theorem \ref{main}.\\

In section \ref{comparison} we prove that our cluster algebra structure agrees with the one given in \cite{gls}. In particular we give the explicit sequence of mutations which relates our initial seed to the one given in \cite{gls}.

\section{Main results}\label{mainresults}
Let $G$ be simple algebraic group of type $G_2$ (over $\C$) and $\Lg$ be its Lie algebra. Let $B^{\pm}$ be a pair of opposite Borel subgroups, $N^{\pm}$ be their respective unipotent radicals. Fix $\{f_1,f_2,h_1,h_2,e_1,e_2\}$ Chevalley generators of $\Lg$. Let  $x_i(t) = \exp{te_i}$ and $y_i(t) = \exp{tf_i}$ (for $i=1,2$) be one parameter subgroups corresponding to Chevalley generators $e_i$ and $f_i$. For each $i=1,2$ let $P_i^-$ be parabolic subgroup generated by $B^-$ and $x_i(t)$. Consider the following multi-homogeneous coordinate rings
\[
\A_i = \C[P_j^-\backslash G]
\]

where $i\in \{1,2\}$ and $j = \{1,2\}\setminus \{i\}$. This notation will be retained throughout these notes.\\

Let $\omega_1, \omega_2$ be fundamental weights of the root system of type $G_2$. It is well known that left $N^-$ invariant regular functions on $G$ is a model for all irreducible finite dimensional representations of $G$, i.e,
\begin{equation}\label{baseaffine}
\C[N^-\backslash G] = \bigoplus_{m,n\in \Z_{\geq 0}} V_{m\omega_1+n\omega_2}
\end{equation}
as right $G$-modules. Here we use the notation $V_{\lambda}$ for (unique) irreducible, finite dimensional $G$ module of highest weight $\lambda$ (where $\lambda$ is some dominant integral weight). With respect to this model, we have
\[
\A_i = \bigoplus_{n\in \Z_{\geq 0}} V_{n\omega_i}
\]

The main result of this note provides a cluster algebra structure on $\A_i$. 

\subsection{Recollections from representation theory}\label{minors}
We have left and right action of $\Lg$ on $\C[G]$. In order to distinguish the two actions, we write $X^{\dagger}$ to denote $X$ acting on right.
\[
(X^{\dagger}f)(x) = \frac{d}{dt}f(x\exp{tX})|_{t=0}
\]
In particular $e_i^{\dagger}$ will be referred to as ``raising operators" and $f_i^{\dagger}$ will be referred to as ``lowering operators". We define certain left $N^-$ invariant functions on $G$, called generalized minors. Following \cite{double}, let $G_0 = N^-HN^+$ be the dense subset of $G$ consisting of elements which admit Gaussian decomposition. For $x\in G_0$ we write $x=[x]_-[x]_0[x]_+$ where $[x]_{\pm}\in N^{\pm}$ and $[x]_0\in H$. For each $i$ define
\[
\Delta^{\omega_i}(x) = [x]_0^{\omega_i}
\]
where $\omega_i$ is thought of as character of torus $H$. Finally one chooses a lifting of $W$ to $G$ by
\[
\overline{s_i} = \varphi_i\left(\begin{array}{cc} 0 & -1 \\ 1 & 0\end{array}\right)
\]
here $\varphi_i : SL_2(\C)\rightarrow G$ is given by:
\[
\varphi_i\left(\begin{array}{cc} 1 & t \\ 0 & 1\end{array}\right) = x_i(t),\ \varphi_i\left(\begin{array}{cc} 1 & 0 \\ t & 1\end{array}\right) = y_i(t)
\]
For any $w\in W$, choose a reduced expression $w = s_{i_1}\ldots s_{i_l}$ and set 
\[
\overline{w} = \overline{s_{i_1}} \ldots \overline{s_{i_l}}
\]
Finally define the function $\Delta_{\omega_i, w\omega_i}$ as
\[
\Delta_{\omega_i,w\omega_i}(x) = \Delta^{\omega_i}(x\overline{w})
\]

Since for most part of this note we will be concerned only with minors of the form $\Delta_{\omega_i, w\omega_i}$ we will denote them by $\Delta^{w\omega_i}$. It is easy to check that each of these functions is left $N^-$ invariant and under the isomorphism (\ref{baseaffine}) we have
\[
\Delta^{w\omega_i} \in V_{\omega_i}(w\omega_i)
\]

Alternatively, in $G_2$ case, one can define $\Delta^{w\omega_i}$ as follows:
\begin{eqnarray}\label{gminors1}
\Delta^{s_1\omega_1} = f_1^{\dagger}\Delta^{\omega_1}\\
\Delta^{s_2s_1\omega_1} = f_2^{\dagger}\Delta^{s_1\omega_1}\\
\Delta^{s_1s_2s_1\omega_1} = \frac{1}{2}(f_1^{\dagger})^2 \Delta^{s_2s_1\omega_1}\\
\Delta^{s_2s_1s_2s_1\omega_1} = f_2^{\dagger}\Delta^{s_1s_2s_1\omega_1}\\
\Delta^{w_0\omega_1} = f_1^{\dagger}\Delta^{s_2s_1s_2s_1\omega_1}
\end{eqnarray}

\begin{eqnarray}\label{gminors2}
\Delta^{s_2\omega_2}  = f_2^{\dagger}\Delta^{\omega_2}\\
\Delta^{s_1s_2\omega_2} = \frac{1}{6}(f_2^{\dagger})^3\Delta^{s_2\omega_2}\\
\Delta^{s_2s_1s_2\omega_2} = \frac{1}{2}(f_2^{\dagger})^2\Delta^{s_1s_2\omega_2}\\
\Delta^{s_1s_2s_1s_2\omega_2} = \frac{1}{6}(f_1^{\dagger})^3\Delta^{s_2s_1s_2\omega_2}\\
\Delta^{w_0\omega_2} = f_2^{\dagger}\Delta^{s_1s_2s_1s_2\omega_2}
\end{eqnarray}

\subsection{Initial seeds for $\A_i$} In this section we define two initial seeds (see Section \ref{review} for terminology)\\

\noindent {\bf Definition of  } $\underline{\Sigma_1}$: Exchange matrix is encoded in following valued graph.

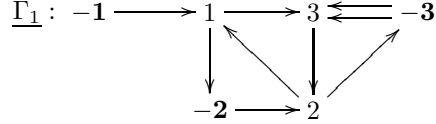
\[
\underline{\Gamma_1}:
\xymatrix{
{\bf -1} \ar[r] & 1 \ar[d] \ar[r] & 3 \ar[d] & {\bf -3} \ar[l]<0.5ex> \ar[l]<-0.5ex> \\
 & {\bf -2} \ar[r] & 2 \ar[ru] \ar[lu] &  
}
\]

The notations of this quiver (and all that we will consider in this paper) are as follows. The bold faced vertices are frozen (correspond to coefficients in the cluster algebra). The arrows between top row and bottom row are valued: if $k$ is a vertex in top row and $l$ is a vertex in bottom row, then $k\rightarrow l$ means: $b_{kl}=3$ and $b_{lk}=-1$. Therefore, in this case, we have following exchange matrix $\underline{B_1}$:

\begin{center}
\begin{tabular}{c|ccc}
& 1 & 2 & 3 \\ \hline
1 & 0 & -3 & 1 \\
2 & 1 & 0 & -1 \\
3 & -1 & 3 & 1 \\
-1 & 1 & 0 & 0\\
-2 & -1 & 1 & 0\\
-3 & 0 & -3 & 2
\end{tabular}
\end{center}

where columns are labeled by $1,2,3$ and rows are labeled by $1,2,3,-1,-2,-3$. \\

In order to  set the functions we define
\[
X_0 = \frac{1}{2} f_1^{\dagger}f_2^{\dagger}f_1^{\dagger} \Delta^{\omega_1}
\]
(this is the zero weight vector in first fundamental representation). Let
\[
X_{-1} = \Delta^{w_0\omega_1} 
\]
\[
X_{-2} = X_0(\Delta^{\omega_1}\Delta^{w_0\omega_1} + \Delta^{s_1\omega_1}\Delta^{s_2s_1s_2s_1\omega_1}) - \Delta^{\omega_1}\Delta^{s_1s_2s_1\omega_1}\Delta^{s_2s_1s_2s_1\omega_1} - \Delta^{s_1\omega_1}\Delta^{s_2s_1\omega_1}\Delta^{w_0\omega_1}\\
\]
\[
X_{-3} = \Delta^{\omega_1} 
\]
\[
X_1 = \Delta^{s_2s_1s_2s_1\omega_1} 
\]
\[
X_2 = \Delta^{s_1\omega_1}(\Delta^{s_2s_1\omega_1})^2 - \Delta^{\omega_1}\Delta^{s_2s_1\omega_1} X_0 - (\Delta^{\omega_1})^2\Delta^{s_2s_1s_2s_1\omega_1} \\
\]
\[
X_3 = \Delta^{s_2s_1\omega_1}
\]

In the terminology of the  theory of cluster algebras, we fix the ground ring to be $\Z[X_{-1}, X_{-2}, X_{-3}^{\pm} ]$ (see section \ref{review} for definitions).\\

\noindent {\bf Definition of } $\underline{\Sigma_2}$: Exchange matrix is encoded in following valued graph
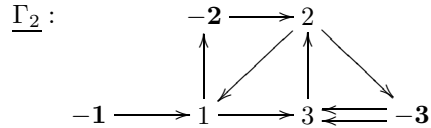
\[ 
\underline{\Gamma_2}:
\xymatrix{
  & {\bf -2} \ar[r] & 2 \ar[ld] \ar[rd] & \\
  {\bf -1} \ar[r] & 1 \ar[u] \ar[r] & 3 \ar[u] & {\bf -3} \ar[l]<0.5ex> \ar[l]<-0.5ex>
}
\]

where the notation remains the same, as in previous part. In this case we have the following exchange matrix, $\underline{B_2}$:

\begin{center}
\begin{tabular}{c|ccc}
 & 1 & 2 & 3 \\ \hline 1&  0 & -1 & 1 \\ 2 & 3 & 0 & -3 \\ 3 & -1 & 1 & 0 \\ -1 & 1 & 0 & 0 \\ -2 & -3 & 1 & 0 \\ -3 & 0 & -1 & 2 
\end{tabular}
\end{center}

where columns are labeled by $\{1,2,3\}$ and rows by $\{1,2,3,-1,-2,-3\}$.Finally we set the functions
\begin{eqnarray*}
Y_{-1} = \Delta^{w_0\omega_2}\\
 Y_{-2} = (2f_1^{\dagger}f_2^{\dagger} - f_2^{\dagger}f_1^{\dagger})Y_2 \\
 Y_{-3} = \Delta^{\omega_2}\\
 Y_1 = \Delta^{s_1s_2s_1s_2\omega_2}\\
 Y_2 = \frac{1}{6}(f_1^{\dagger})^2f_2^{\dagger}\Delta^{\omega_2}\\
 Y_3 = \Delta^{s_1s_2\omega_2}
\end{eqnarray*}

Again we fix ground ring to be $\Z[Y_{-1},Y_{-2},Y_{-3}^{\pm}]$.
\begin{rem}\begin{rm}
The choice of the functions corresponding to cluster variables will be explained in section \ref{explanation}.
\end{rm}\end{rem}

\begin{rem}\begin{rm}
We will see in section \ref{comparison} that our initial seed(s) and the ones given in \cite{gls} give rise to isomorphic cluster algebras. This fact combined with Theorem \ref{main} will give an affirmative answer to Conjecture 10.4 of \cite{gls}.
\end{rm}\end{rem}

\subsection{Statements of main results}
\begin{thm}\label{main}
The cluster algebra associated to initial seed $\underline{\Sigma_i}$ is isomorphic to localization of $\A_i$ at multiplicative subset $S_i = \{(\Delta^{\omega_i})^n\}_{n\geq 0}$.
\[
\A(\underline{\Sigma_i})_{\C} = S_i^{-1}\A_i
\]
\end{thm}

We note few corollaries to this theorem.

\begin{cor}\label{corollary1}
The upper cluster algebra structure on $\C[G^{e,w_0}]$ (as given in \cite{CA3}, see Section \ref{review} for definition) is in fact a cluster algebra structure.
\end{cor}

As remarked in \cite{gls} (see discussion following Conjecture 10.4), affirmative statement of Theorem \ref{main} will imply the following (Conjecture 9.6 of \cite{gls}):

\begin{cor}\label{corollary2}
The projection $ \A_i \rightarrow \C[N_i]$ induced from inclusion $N_i \rightarrow P_j^-\backslash G$ (where $\{j\} = \{1,2\}\setminus \{i\}$)  defines cluster algebra structure on $\C[N_i]$
\end{cor}

\section{Recollections from \cite{CA3}}\label{review}
In this section we will recall some results from \cite{CA3}. We will mainly focus on description of initial seed for the coordinate rings of double Bruhat cells, in terms of double reduced words.

\subsection{Definitions from theory of cluster algebras}
The definitions we recall here are not the most general ones. For the general theory of cluster algebras, reader is referred to papers \cite{CA1}, \cite{CA3} and \cite{CA4}.\\

Let $m \geq n$ be two non-negative integers. Let $\F = \Q(u_1,\ldots, u_m)$ be field of rational functions in $m$ variables. A {\em seed} in $\F$ is a tuple $\Sigma = (\widetilde{{\bf x}}, \widetilde{B})$ where
\begin{itemize}
\item[-] $\widetilde{{\bf x}} = \{x_1,\ldots, x_m\}$ is a set of $m$ free generators  of $\F$ over $\Q$.
\item[-] $\widetilde{B} = (b_{ij})$ is $m\times n$ matrix with integer entries, such that the square submatrix $B = (b_{ij})_{1\leq i,j\leq n}$ obtained by taking first $n$ rows (called principal part of $\widetilde{B}$) is skew symmetrizable (i.e, there exists a diagonal matrix $D$ with positive diagonal entries, such that $DB$ is skew symmetric).
\end{itemize}

For each $k$, $1\leq k\leq n$, we define {\em mutation of $\Sigma$ in direction $k$} as another seed $\Sigma' = \mu_k(\Sigma)$ if $\Sigma' = (\widetilde{{\bf x'}}, \widetilde{B}')$ where
\begin{itemize}
\item[-] $\widetilde{{\bf x'}} = (x_1',\ldots, x_m')$ are given by
\[
x_i' = x_i \mbox{ if } i\not= k
\]
\[
x_kx_k' = \Pi_{i=1}^m x_i^{[b_{ik}]_+} + \Pi_{i=1}^m x_i^{[-b_{ik}]_+}
\]
here we use the notation $[b]_+ = \mbox{max}(0,b)$.
\item[-] The matrix $\widetilde{B}' = (b_{ij}')$ is given by
\[
b_{ij}' = \left\{\begin{array}{lcr} -b_{ij} & \mbox{if} & i=k \mbox{ or } j=k \\ b_{ij} + \mbox{sgn}(b_{ik})[b_{ik}b_{kj}]_+ & \mbox{otherwise} & \end{array}\right.
\]
\end{itemize}

It is well known (and easy to check) that a) the principal part of $\widetilde{B}'$ is again skew-symmetrizable; and b) mutation is an involution (i.e, $\mu_k(\mu_k(\Sigma)) = \Sigma$). This allows us to define equivalence relation: two seeds $\Sigma$ and $\Sigma'$ in $\F$ are {\em mutation equivalent} if there exists a sequence $(k_1,\ldots, k_l)$ such that $\Sigma' = \mu_{k_1}\ldots \mu_{k_l}(\Sigma)$. \\

If $\Sigma = (\widetilde{{\bf x}},\widetilde{B})$ is a seed of $\F$, then the set $\widetilde{{\bf x}}$ is called {\em extended cluster}, ${\bf x} = \{x_1,\ldots, x_n\}$ is called cluster and $x_i$, $1\leq i\leq n$ are called {\em cluster variables}. \\

Note that by applying mutations we only change the functions $x_i$ $1\leq i\leq n$. Therefore if $\Sigma$ is mutation equivalent to $\Sigma'$, then $x_j = x_j'$ for every $n<j\leq m$. We refer to $x_j$ with $n<j\leq m$ as {\em coefficients}. We fix the ground ring $R$ as some subring of $\Z[x_j^{\pm} : n<j\leq m]$ containing $\Z[x_j : n<j\leq m]$, i.e,
\[
\Z[x_{n+1},\ldots, x_m] \subseteq R \subseteq \Z[x_{n+1}^{\pm},\ldots, x_m^{\pm}]
\]

Define {\em the upper cluster algebra} $\overline{\A(\Sigma)}$ as subring of $\F$ consisting of functions $f\in \F$ such that for every $\Sigma' \sim \Sigma$, $f$ can be expressed as Laurent polynomial in cluster variables of $\Sigma'$ with coefficients from $R$.\\

We define cluster algebra associated to seed $\Sigma$, denoted by $\A(\Sigma)$ as $R$-subalgebra of $\F$ generated by all cluster variables belonging to seeds $\Sigma' \sim \Sigma$. The following result is known as {\em Laurent phenomenon}
\[
\A(\Sigma) \subset \overline{\A(\Sigma)}
\]

We say that a seed $\Sigma = (\widetilde{{\bf x}},\widetilde{B})$ is acyclic if the directed graph with vertex set $\{1,\ldots, n\}$ and arrow $k\rightarrow l$ if and only if $b_{kl}>0$, has no oriented cycles. It is known that if in the mutation equivalence class of $\Sigma$, there is some acyclic seed and the matrix $\widetilde{B}$ has full rank, then $\A(\Sigma) = \overline{\A(\Sigma)}$.\\

Similarly we say that seed $\Sigma$ is bipartite if the directed graph defined by $B$ as above is bipartite. \\

\subsection{Double Bruhat cells and cluster algebras} Given a pair of elements $u,v\in W$, the double Bruhat cell $G^{u,v}$ is defined as
\[
G^{u,v} = BuB \cap B^-vB^-
\]

In \cite{CA3} authors give a combinatorial way to construct seed $\Sigma$ of $\F = \C(G^{u,v})$, starting from a double reduced word of pair $(u,v)\in W\times W$. We will state the results of \cite{CA3} regarding coordinate rings of double Bruhat cells, only in our context.\\

We again let $G$ to be simple algebraic group over $\C$ of type $G_2$. Let $w_0$ be longest element of Weyl group of $G$. The relevant double Bruhat cells is
\[
G^{e,w_0} = B\cap B^-w_0B^-
\]

As application of Proposition 2.8 of \cite{CA3} we have

\begin{prop}
The subset $G^{e,w_0}$ is defined in $G$ by following (recall definitions of generalized minors from section \ref{minors}):
\[
\Delta_{u\omega_i, \omega_i} = 0,\ u\not=  e
\]
\[
\Delta_{\omega_i,\omega_i} \not= 0,\ \ \Delta_{\omega_i, w_0\omega_i} \not= 0
\]
\end{prop}

The inclusion $G^{e,w_0} \rightarrow G$ composed with projection onto $N^-\backslash G$ allows us to consider $\C[N^-\backslash G]$ as subalgebra of $\C[G^{e,w_0}]$. In other words, $\C[G^{e,w_0}]$ is isomorphic to localization of $\C[N^-\backslash G]$ at functions $\Delta^{\omega_i}$ and $\Delta^{w_0\omega_i}$. \\

We have two reduced words for $w_0$, namely ${\bf i_1} = (1,2,1,2,1,2)$ and ${\bf i_2} = (2,1,2,1,2,1)$. Each one of them defines an initial seed in $\C(G^{e,w_0})$ (see section 2 of \cite{CA3}):\\

\noindent {\bf Definition of } $\Sigma_1$: The exchange matrix $\widetilde{B}_1$ is encoded in following valued quiver
\[
\Gamma_1:
\xymatrix{
{\bf -1} \ar[r] & 1 \ar[r] \ar[ld]& 3 \ar[r] \ar[ld]& {\bf 5}\ar[ld] \\
{\bf -2} \ar[r] & 2 \ar[r] \ar[u]& 4 \ar[r] \ar[u]& {\bf 6}
}
\]

The functions corresponding to cluster variables are:
\[
X_{-1} = \Delta^{w_0\omega_1},\ X_1 = \Delta^{s_2s_1s_2s_1\omega_1},\ X_3 = \Delta^{s_2s_1\omega_1},\ X_5 = \Delta^{\omega_1}
\]
\[
X_{-2} = \Delta^{w_0\omega_2},\ X_2 = \Delta^{s_2s_1s_2\omega_2},\  X_4 = \Delta^{s_2\omega_2},\  X_6 = \Delta^{\omega_2}
\]

\noindent {\bf Definition of } $\Sigma_2$: The exchange matrix $\widetilde{B}_2$ is encoded by following valued quiver
\[
\Gamma_2:
\xymatrix{
{\bf -1} \ar[r] & 2 \ar[r] \ar[d]& 4 \ar[r] \ar[d]& {\bf 6} \\
{\bf -2} \ar[r] & 1 \ar[r] \ar[lu]& 3 \ar[r] \ar[lu]& {\bf 5} \ar[lu]
}
\]

The functions corresponding to cluster variables are:
\[
X_{-1} = \Delta^{w_0\omega_1},\ X_2 = \Delta^{s_1s_2s_1\omega_1},\ X_4 = \Delta^{s_1\omega_1},\ X_6 = \Delta^{\omega_1}
\]
\[
X_{-2} = \Delta^{w_0\omega_2},\ X_1 = \Delta^{s_1s_2s_1s_2\omega_2},\ X_3 = \Delta^{s_1s_2\omega_2},\ X_5 = \Delta^{\omega_2}
\]

We have the following special case of Theorem 2.10 of \cite{CA3}:

\begin{thm}\label{bfz}
The upper cluster algebra $\overline{\A(\Sigma_i)}_{\C}$ is isomorphic to $\C[G^{e,w_0}]$. In particular every function obtained by applying sequence of mutations to $\Sigma_i$ is regular function on $G^{e,w_0}$.
\end{thm}

\begin{rem}\begin{rm}
It can be shown that the upper cluster algebras $\overline{\A(\Sigma_i)}$ are in fact {\em equal}, i.e, $\Sigma_1$ is mutation equivalent to $\Sigma_2$ by following sequence of mutations:
\[
\Sigma_1 = \mu_1\mu_3\mu_2\mu_1\mu_2\mu_4\mu_2\mu_1\mu_2\mu_3\mu_1(\Sigma_2)
\]
In order to get equality above, one needs to relabel cluster variables after the sequence of mutations. This can be verified by using quiver mutation software of Prof. B. Keller and verifying the determinantal identities for $G_2$ given in \cite{schubert}. This result will not be used in this paper.
\end{rm}\end{rem}

\begin{rem}\begin{rm}
It will follow from our main results that we have equality
\[
\A(\Sigma_i) = \overline{\A(\Sigma_i)}
\]
\end{rm}\end{rem}

\section{Explanation of choice of seeds for $\A_i$}\label{explanation}
In this section we explain the ``greedy approach" to obtain initial seeds of $\A_i$, starting from the initial seeds of section \ref{review}. The idea is to apply minimum number of mutations to get the cluster variables belonging to single $\A_i$.\\

\noindent $\A_1$: We start from initial seed $\Sigma_1$ and apply mutations $\mu_2\mu_4$ to obtain the following:
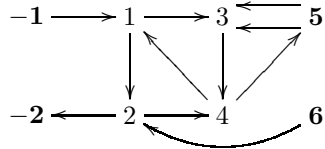
\[
\xymatrix{
{\bf -1} \ar[r] & 1 \ar[r] \ar[d]& 3 \ar[d] & {\bf 5} \ar[l]<1ex> \ar[l]<-1ex> \\
{\bf -2}  & 2 \ar[l] \ar[r] & 4 \ar[ru] \ar[lu]& {\bf 6} \ar @/^1pc/[ll]
}
\]

\begin{lem}\label{seed1}
Let $X_2'$ and $X_4'$ be functions obtained by applying mutations $\mu_2\mu_4$ to initial seed $\Sigma_1$. Then we have
\[
X_2' = X_0(\Delta^{\omega_1}\Delta^{w_0\omega_1}+\Delta^{s_1\omega_1}\Delta^{s_2s_1s_2s_1\omega_1}) - \Delta^{\omega_1}\Delta^{s_1s_2s_1\omega_1}\Delta^{s_2s_1s_2s_1\omega_1} - \Delta^{s_1\omega_1}\Delta^{s_2s_1\omega_1}\Delta^{w_0\omega_1}
\]
\[
X_4' = \Delta^{s_1\omega_1}(\Delta^{s_2s_1\omega_1})^2 - \Delta^{\omega_1}\Delta^{s_2s_1\omega_1}X_0 - (\Delta^{\omega_1})^2\Delta^{s_2s_1s_2s_1\omega_1}
\]
\end{lem}

This lemma proves that the initial seed $\underline{\Sigma_1}$ is precisely the one obtained from $\mu_2\mu_4(\Sigma_1)$ by ``freezing" the vertex labeled $2$ (and renaming vertices, just for convenience of notation). In the mutated seed $\mu_2\mu_4(\Sigma_1)$ all cluster variables belong to representations of type $V_{n\omega_1}$ except for ones corresponding to vertices ${\bf -2}$ and ${\bf 6}$ which are only linked to vertex $2$. Therefore if we do not allow mutation at this vertex, all the functions we shall obtain will again belong to $V_{n\omega_1}$ (for some $n$) and hence to $\A_1$.\\

\noindent $\A_2$: In this case we start with initial seed $\Sigma_2$ and apply mutations $\mu_2\mu_4$ to obtain:
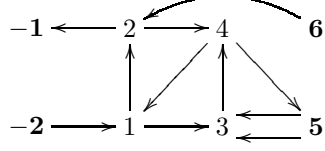
\[
\xymatrix{
{\bf -1} & 2 \ar[l] \ar[r] & 4 \ar[ld] \ar[rd] & {\bf 6} \ar @/_1pc/[ll] \\
{\bf -2} \ar[r] & 1 \ar[u] \ar[r] & 3 \ar[u] & {\bf 5} \ar[l]<-1ex> \ar[l]<1ex>
}
\]

\begin{lem}\label{seed2}
Let $X_2'$ and $X_4'$ be functions obtained by applying mutations $\mu_2\mu_4$ to initial seed $\Sigma_2$. Then we have
\[
X_2' = (2f_1^{\dagger}f_2^{\dagger} - f_2^{\dagger}f_1^{\dagger})X_4'
\]
\[
X_4' = \frac{1}{6} (f_1^{\dagger})^2f_2^{\dagger}\Delta^{\omega_2}
\]
\end{lem}

We will prove these lemmas in sections \ref{lemmas}

\section{Proof of Theorem \ref{main}}\label{proof}

In the next two subsections we describe the main constituents of the proof of Theorem \ref{main}. We treat the cases $i=1$ and $i=2$ separately.

\subsection{Proof for case $i=1$}\label{proof1}

We first assume that ground ring for cluster algebra $\A(\underline{\Sigma_1})$ is $\Z[X_{-1}, X_{-2}, X_{-3}]$. \\

Lemma \ref{seed1} implies  that cluster algebra $\A(\underline{\Sigma_1})$ is a subalgebra of $\A(\Sigma_1)$, that is the functions obtained (as cluster variables) by applying mutations to $\underline{\Sigma_1}$ are cluster variables of $\A(\Sigma_1)$. This allows us to apply Theorem \ref{bfz} to conclude that

\[
\A(\underline{\Sigma_1}) \subset \A_1
\]

We define the following $\Z_{\geq 0}\times \Z\times \Z$ grading on $\A(\Sigma_1)$. An element $F\in \A(\underline{\Sigma_1})$ has degree $(n,p,q)$ if considered as an element of $\A_1$, we have
\[
F\in V_{n\omega_1}(p\alpha_1 + q\alpha_2)
\]

With this definition we have degrees of cluster variables of initial seed $\underline{\Sigma_1}$:
\[
\mbox{deg}(X_{-1}) =  (1,-2,-1),\ \mbox{deg}(X_{-2}) = (3,0,0),\ \mbox{deg}(X_{-3}) = (1,2,1)
\]
\[
\mbox{deg}(X_1) = (1,-1,-1),\ \mbox{deg}(X_2) = (3,3,1),\  \mbox{deg}(X_3) = (1,1,0)
\]

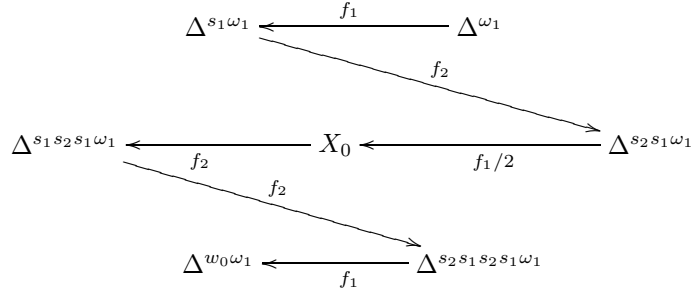
\begin{figure}\caption{First fundamental representation}\label{first}
\[
\xymatrix@R=0.433in@C=0.25in{
& \Delta^{s_1\omega_1} \ar[rrrd]^{f_2}& & \Delta^{\omega_1} \ar[ll]_{f_1} & \\
\Delta^{s_1s_2s_1\omega_1} \ar[rrrd]^{f_2}& & X_0 \ar[ll] ^{f_2}& & \Delta^{s_2s_1\omega_1} \ar[ll]^{f_1/2} \\
& \Delta^{w_0\omega_1} & & \Delta^{s_2s_1s_2s_1\omega_1}\ar[ll]^{f_1} &
}
\]
\end{figure}

\noindent Recall that the first fundamental representation $V_{\omega_1}$ has the realization in terms of generalized minors (see Figure \ref{first})\\

Since the algebra $\A_1$ is generated by the weight vectors of $V_{\omega_1}$, in order to prove the first part of Theorem \ref{main}, we need to prove that the inclusion $\A(\underline{\Sigma_1})\subset \A_1$ becomes equality when localized at $\Delta^{\omega_1}$. Thus it would suffice to prove the following two statements.

\begin{itemize}
\item[a)] There are cluster variables $W, Y$ such that
\[
\mbox{deg}(W) = (1,-1,0),\ \mbox{deg}(Y) = (1,1,1)
\]
\item[b)] $X_0\Delta^{\omega_1}$ can be written in terms of cluster variables (i.e, it belongs to algebra generated by cluster variables).
\end{itemize}

We begin by applying $\mu_3$ to $\underline{\Sigma_1}$. Let $X_3'$ be cluster variable obtained at vertex $3$.
\[
X_3' = \frac{X_1X_{-3}^2+X_2}{X_3} = \Delta^{s_1\omega_1}\Delta^{s_2s_1\omega_1} - \Delta^{\omega_1}X_0
\]

Therefore we have $\Delta^{\omega_1}X_0$ in terms of other cluster variables (assuming that $\Delta^{s_1\omega_1}$ and $\Delta^{s_2s_1\omega_1}$ are cluster variables: the assertion of part a)).\\

Let us define $\underline{\Sigma_1}^0 := \mu_2\mu_3(\underline{\Sigma_1})$. The exchange matrix at this cluster is given by

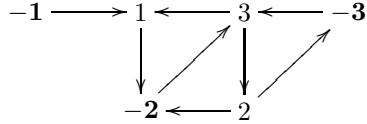
\[
\xymatrix{
{\bf -1} \ar[r] & 1 \ar[d] & 3 \ar[l] \ar[d] & {\bf -3} \ar[l] \\
& {\bf -2} \ar[ru] & 2 \ar[l] \ar[ru] &  
}
\]

Therefore the seed $\underline{\Sigma_1}^0$ is bipartite. Following \cite{CA4} we define
\[
\mu_+ = \mu_1\mu_2,\ \ \mu_- = \mu_3
\]
As long as we restrict ourselves to ``the bipartite belt" the mutations $\mu_1$ and $\mu_2$ commute (since vertices $1$ and $2$ are not linked). Further we define
\[
\underline{\Sigma_1}^r = \left\{\begin{array}{rcl} \underbrace{\cdots \mu_-\mu_+\mu_-\mu_+}_{r\mbox{ terms}} (\underline{\Sigma_1}^0) & if & r\geq 0 \\
\underbrace{\cdots\mu_+\mu_-\mu_+\mu_-}_{-r\mbox{ terms}}(\underline{\Sigma_1}^0) & if & r<0 
\end{array}\right.
\]

We denote by $X_i^{(r)}$, the cluster variables at seed $\underline{\Sigma_1}^r$ and $d_i^{(r)} \in \N\times \Z\times \Z$ its degree. Then part a) follows from following proposition: to see that there exist cluster variables with degrees $(1,-1,0)$ and $(1,1,1)$ we just observe that $d_1^{(-3)} = (1,1,1)$ and $d_1^{(-7)} = (1,-1,0)$. This proposition will be proved in section \ref{propositions}.

\begin{prop}\label{degrees1}
In the notation introduced above, we have
\[
d_1^{(2r)} = d_1^{(2r-1)},\ \ d_2^{(2r)} = d_2^{(2r-1)},\ \ d_3^{(2r+1)} = d_3^{(2r)}
\]
\[
d_1^{(2r+1)} = \left\{\begin{array}{lcl}
(r/2+1, 1,1) & if & r\geq 0 \mbox{ and even} \\
((r+1)/2, 1,0) & if & r>0 \mbox{ and odd} \\
(1,1,1) & if & r=-2 \\
((1-r)/2,-1,-1) & if & r<0 \mbox{ and odd}\\
(-r/2-1, -1, 0) & if & r<-2 \mbox{ and even}
\end{array}\right.
\]
\[
d_2^{(2r+1)} = \left\{ \begin{array}{lcl} 
(3r/2+3, 3,1) & if & r\geq0 \mbox{ and even }\\ 
(3(r+1)/2+3,3,2) & if & r>0 \mbox{ and odd} \\
(3,3,2) & if & r=-1 \\
(-3r/2, -3, -2) & if & r<0 \mbox{ and even}\\
(-3(r+1)/2, -3,-1) & if & r<-1 \mbox{ and odd}
\end{array}\right.
\]
\[
d_3^{(2r)} = \left\{\begin{array}{lcl} (2+r, 2,1) & \mbox{if} & r\geq 0 \\ (2,0,0) & \mbox{if} & r=-1 \\ (-r, -2,-1) & \mbox{if} & r<-1 \end{array}\right.
\]
\end{prop}

\newpage

\subsection{Proof for the case $i=2$}\label{proof2}

The proof of this part is exactly similar to the previous one. We restrict ourselves to ground ring $\Z[Y_{-1}, Y_{-2}, Y_{-3}]$. Assuming Lemma \ref{seed2} we conclude from Theorem \ref{bfz} the inclusion
\[
\A(\underline{\Sigma_2}) \subset \A_2
\]

which allows us to define a grading on elements of $\A(\underline{\Sigma_2})$. We say $\mbox{deg}(F) = (n,p,q)$ if as an element of $\A_2$ we have
\[
F \in V_{n\omega_2}(p\alpha_1+q\alpha_2)
\]

In this notation the cluster variables have following degrees
\[
\mbox{deg}(Y_{-1}) = (1,-3,-2),\ \ \mbox{deg}(Y_{-2}) = (1,0,0),\ \  \mbox{deg}(Y_{-3}) = (1,3,2) 
\]
\[
\mbox{deg}(Y_1) = (1,-3,-1),\ \ \mbox{deg}(Y_2) = (1,1,1),\ \ \mbox{deg}(Y_3) = (1,0,1)
\]

\begin{figure}\caption{Second Fundamental Representation}\label{second}
\[
\xymatrix@R=0.3in@C=0.0625in{
&&& \Delta^{\omega_2} \ar[rrrd]^{f_2}&&& \\
\Delta^{s_1s_2\omega_2}\ar[rrrd]^{f_2 = (1\ 2)^T} && F(1,1) \ar[ll]_{f_1} \ar[rrrd]^{f_2}&& F(2,1) \ar[ll]_{f_1/2}&& \Delta^{s_2\omega_2}\ar[ll]_{f_1/3} \\
& F(-1,0)\ar[rrrd]^{f_2} & & *\txt{$F^1(0,0)$\\ $F^2(0,0)$} \ar[ll]_{f_1=(1\ 1)} \ar[rrrd]^{f_2=(0\ 1)}&& F(1,0)\ar[ll]_{f_1 = (1\ 1)^T} & \\
\Delta^{s_1s_2s_1s_2\omega_2} \ar[rrrd]_{f_2}&& F(-2,-1)\ar[ll]^{f_1} && F(-1,-1)\ar[ll]^{f_1/2} && \Delta^{s_2s_1s_2\omega_2} \ar[ll]^{f_1/3}\\
&&& \Delta^{w_0\omega_2} &&&
}
\]
\end{figure}
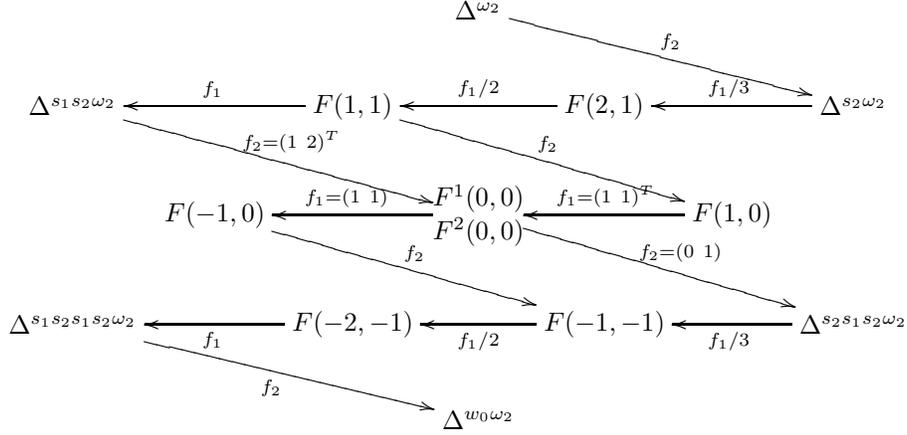

\noindent Note that the second fundamental representation of $G$ has the  realization in terms of certain functions on $G$ (see Figure \ref{second}). 

\noindent In the figure above, each function represents a chosen basis vector of the respective weight space and the arrows indicate action of lowering operators. For instance $\xymatrix{ \Delta^{s_1s_2\omega_2} \ar[rr]^{f_2 = (1\ 2)^T} & & *\txt{$F^1(0,0)$\\$F^2(0,0)$} }$ means that $f_2^{\dagger}(\Delta^{s_1s_2\omega_2}) = F^1(0,0) + 2F^2(0,0)$. The functions $F(i,j)$ have weight $\alpha_1 i+ \alpha_2 j$ and Figure \ref{second} can be taken as definition of these functions. For example

\[
F(1,1) = \frac{1}{6}(f_1^{\dagger})^2 \Delta^{s_2\omega_2} = Y_2
\]
\[
F^1(0,0) = (2f_1^{\dagger}f_2^{\dagger} - f_2^{\dagger}f_1^{\dagger})(F(1,1)) = Y_{-2}
\]

Again it suffices to prove the following two statements
\begin{itemize}
\item[a)] There exist cluster variables with degrees $(1,3,1); (1,2,1); (1,1,0); (1,-1,0); (1,0,-1); (1,-1,-1)$ and $(1,-2,-1)$.
\item[b)] $\Delta^{\omega_2}.F^2(0,0)$ can be written as a polynomial in cluster variables.
\end{itemize}

Similar to previous part, we define $\underline{\Sigma_2}^0$ to be $\mu_2\mu_3(\underline{\Sigma_2})$ to make it bipartite:
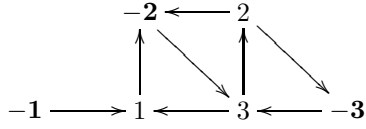
\[
\xymatrix{
 & {\bf -2} \ar[rd] & 2 \ar[l] \ar[rd] & \\
 {\bf -1} \ar[r] & 1 \ar[u] & 3 \ar[l] \ar[u] & {\bf -3} \ar[l]
}
\]

And define 
\[
\mu_+ = \mu_1\mu_2,\ \ \mu_- = \mu_3
\]

The bipartite belt consists of following seeds

\[
\underline{\Sigma_2}^r = \left\{\begin{array}{rcl} \underbrace{\cdots \mu_-\mu_+\mu_-\mu_+}_{r\mbox{ terms}} (\underline{\Sigma_2}^0) & if & r\geq 0 \\
\underbrace{\cdots\mu_+\mu_-\mu_+\mu_-}_{-r\mbox{ terms}}(\underline{\Sigma_2}^0) & if & r<0 
\end{array}\right.
\]

Again we let $Y_i^{(r)}$ be the cluster variables in $\underline{\Sigma_2}^0$ and let $g_i^{(r)}$ be its degree. Also let $U$ and $Z$ be cluster variables appearing at vertex $2$ in $\mu_2\mu_3\mu_1(\underline{\Sigma_2}^0)$ and $\mu_2\mu_3\mu_2(\underline{\Sigma_2}^0)$. Then part a) follows from following degree computation: again we only need to observe that $g_1^{(-3)} = (1,3,1)$, $g_2^{(-3)} = (1,2,1)$, $g_1^{(-7)}=(1,0,-1)$, $g_2^{(-5)} = (1,-1,-1)$, $g_2^{(-3)} = (1,-2,-1)$. This together with the first statement of the proposition exhausts the list of weights demanded in part a). Again the proof of this proposition is given in section \ref{propositions}.

\begin{prop}\label{degrees2}
We have
\[
\mbox{deg}(U) = (1,1,0),\ \ \mbox{deg}(Z) = (1,-1,0)
\]
Moreover we have
\[
g_1^{(2r)} = g_1^{(2r-1)},\ \ g_2^{(2r)} = g_2^{(2r-1)},\ \ g_3^{(2r+1)} = g_3^{(2r)}
\]

\[
g_1^{(2r+1)} = \left\{\begin{array}{lcl}
(r/2+2, 3,1) & \mbox{if} & r\geq 0 \mbox{ and even} \\
((r+1)/2, 0,1) & \mbox{if} & r>0 \mbox{ and odd}\\
(1,3,1) & \mbox{if} & r=-2\\
((1-r)/2, -3,-1) & \mbox{if} & r<0 \mbox{ and odd}\\
(-r/2-1, 0,-1) & \mbox{if} & r<-2 \mbox{ and even}
\end{array}\right.
\]

\[
g_2^{(2r+1)} = \left\{\begin{array}{lcl}
(r/2+1, 2,1) & \mbox{if} & r\geq0 \mbox{ and even}\\
((r+3)/2, 1,1) & \mbox{if} & r>0 \mbox{ and odd}\\
(1,2,1) & \mbox{if} & r=-1 \\
(-r/2,-2,-1) & \mbox{if} & r<0 \mbox{ and even}\\
(-(r+1)/2, -1,-1) & \mbox{if} & r<-1 \mbox{ and odd}
\end{array}\right.
\]
\[
g_3^{(2r)} = \left\{\begin{array}{lcl}
(2+r, 3,2) & \mbox{if} & r \geq 0\\
(2,0,0) & \mbox{if} & r=-1\\
(-r, -3,-2) & \mbox{if} & r<-1 
\end{array}\right.
\]
\end{prop}

Part b) will follow from above proposition together with following ``Pl\"{u}cker relation" (see Section \ref{plucker}):
\[
\Delta^{\omega_2}(F^1(0,0)+F^2(0,0)) = \Delta^{s_2\omega_2}\Delta^{s_1s_2\omega_2} - F(1,1)F(2,1)
\]

\begin{rem}\begin{rm}
Thus to prove the main result we are reduced to proving two lemmas (\ref{seed1} and \ref{seed2}) and two propositions (\ref{degrees1} and \ref{degrees2}). Both the propositions are very combinatorial in nature and are proved in next section. It turns out that degrees $d_i^{(r)}$ and $g_i^{(r)}$ satisfy the same recurrence relations and just have different initial values, which allows us to prove them simultaneously. The lemmas stumble upon certain ``determinant identities" which we prove using representation theory of $G_2$ in Section \ref{lemmas}.
\end{rm}\end{rem}

\begin{rem}\begin{rm}
It follows from the general theory of cluster algebras of affine type (the reader can consult \cite{giovanni} for some structural results) that $U$, $Z$ and $Y_i^{(r)}$ are all the cluster variables  of cluster algebra $\A(\underline{\Sigma_2})$ (and similarly for $\A(\underline{\Sigma_1})$).
\end{rm}\end{rem}

\section{Proofs of Propositions \ref{degrees1} and \ref{degrees2}}\label{propositions}
This section is devoted to the computation of degrees. In order to state the problem at hand in purely combinatorial terms, we review the set up a little bit.\\

Consider the following quivers (transpose of each other):
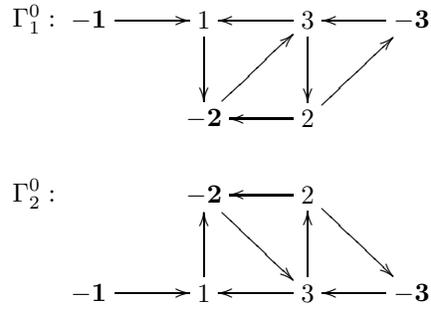
\[
\Gamma_1^0 : 
\xymatrix{
{\bf -1} \ar[r] & 1 \ar[d] & 3 \ar[l] \ar[d] & {\bf -3} \ar[l] \\
 & {\bf -2} \ar[ru] & 2 \ar[l] \ar[ru] & 
}
\]
\[
\Gamma_2^0 : 
\xymatrix{
 & {\bf -2} \ar[rd] & 2 \ar[l] \ar[rd] & \\
 {\bf -1} \ar[r] & 1 \ar[u] & 3 \ar[l] \ar[u] & {\bf -3} \ar[l]
}
\]

We apply the mutation at vertices $1, 2$ or $3$. Let us call arrows between vertices of same parity as {\em horizontal} and between different parity {\em vertical}. The following list describes the mutation rules for $\mu_i$ (note the appearance of factor $3$ to accommodate the valuations on vertical arrows) 

\begin{itemize}
\item[(M1)] If $j\rightarrow i\rightarrow k$ are not both vertical, then add one arrow $j\rightarrow k$.
\item[(M1')] If $j\rightarrow i\rightarrow k$ are both vertical then add three arrows from $j$ to $k$.
\item[(M2)] If $i\rightarrow j$ is an arrow (horizontal or vertical) then change it to $i\leftarrow j$ (similarly other way around).
\item[(M3)] Finally remove all arrows between vertices $-1, -2$ and $-3$ and cancel two cycles, i.e, if there are $s$ arrows from $j$ to $k$ and $t$ arrows from $k$ to $j$ and $s\geq t$ then just write $s-t$ arrows from $j$ to $k$. 
\end{itemize}

It is clear that if two vertices $i$ and $j$ are not connected then $\mu_i\mu_j = \mu_j\mu_i$. This allows us to unambiguously define $\mu_+ = \mu_1\mu_2$ and $\mu_- = \mu_3$. Set
\[
\Gamma_i^r := \left\{\begin{array}{lcl} 
\underbrace{\cdots\mu_-\mu_+\mu_-\mu_+}_{r\mbox{ terms}} \Gamma_i^0 & \mbox{if} & r\geq 0\\
\underbrace{\cdots\mu_+\mu_-\mu_+\mu_-}_{-r\mbox{ terms}} \Gamma_i^0 & \mbox{if} & r<0
\end{array}\right.
\]

\subsection{Structure of graphs $\Gamma_i^r$}
It is clear that between vertices $1,2$ and $3$ the above graphs have following structure:
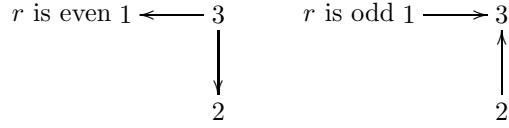
\[
r \mbox{ is even}
\xymatrix{
1 & 3 \ar[l] \ar[d] \\
& 2
}
\ \ \ \ \ \ \ \  
r\mbox{ is odd}
\xymatrix{
1\ar[r] & 3 \\
& 2\ar[u]
}
\]

Since there are no arrows between vertices $-1,-2$ and $-3$, in order to completely determine the structure of graphs $\Gamma_i^r$ we need to compute the following matrix
\[
C_{ij}^{(r)} := \mbox{ Number of arrows from $i$ to $j$}
\]
where $i\in \{-1,-2,-3\}$ and $j\in \{1,2,3\}$. It is clear that both the graphs $\Gamma_1$ and $\Gamma_2$ have same $C$-matrix.\\

The mutation rules define following recurrence relations among entries of $C^{(r)}$ (where the notation $[x]_- := [-x]_+ = \mbox{max}(0,-x)$ is used):

\[
C^{(r+1)} = \left\{\begin{array}{lcl}
\left( \begin{array}{ccc} 
-C_{-1,1}^{(r)} & -C_{-1,2}^{(r)} & C_{-1,3}^{(r)}-[C_{-1,1}^{(r)}]_- - 3[C_{-1,2}^{(r)}]_- \\
-C_{-2,1}^{(r)} & -C_{-2,2}^{(r)} & C_{-2,3}^{(r)}-[C_{-2,1}^{(r)}]_- - [C_{-2,2}^{(r)}]_-\\
-C_{-3,1}^{(r)} & -C_{-3,2}^{(r)} & C_{-3,3}^{(r)}-[C_{-3,1}^{(r)}]_- - 3[C_{-3,2}^{(r)}]_-
\end{array}\right) & \mbox{if} & r \mbox{ is even} \\
 & & \\
\left( \begin{array}{ccc}
C_{-1,1}^{(r)}-[C_{-1,3}^{(r)}]_- & C_{-1,2}^{(r)} - [C_{-1,3}^{(r)}]_- & -C_{-1,3}^{(r)} \\
C_{-2,1}^{(r)}-[C_{-2,3}^{(r)}]_- & C_{-2,2}^{(r)} - 3[C_{-2,3}^{(r)}]_- & -C_{-2,3}^{(r)}\\
C_{-3,1}^{(r)}-[C_{-3,3}^{(r)}]_- & C_{-3,2}^{(r)} - [C_{-3,3}^{(r)}]_- & -C_{-3,3}^{(r)}
\end{array} \right) & \mbox{if} & r \mbox{ is odd}
\end{array}\right.
\]
with the initial value
\[
C^{(0)} = \left( \begin{array}{ccc} 1 & 0 & 0 \\ -1 & -1 & 1 \\ 0 & -1 & 1\end{array}\right)
\]

One can verify directly that the following lemma gives the solution of this recurrence relation:

\begin{lem}
Let $C_i^{(r)}$ be the row of $C$ labeled by $i$. Then we have following solution
\[
C_{-1}^{(r)} = -C_{-1}^{(-r+1)}
\]
\[
C_{-1}^{(r)} = \left(C_{-1,1}^{(r)}, (-1)^{r+1} \left\lfloor \frac{r}{4} \right\rfloor,(-1)^r \left(\left\lceil \frac{r}{2}\right\rceil -1 \right)  \right)\ \mbox{ for } r\geq 1
\]
\[
C_{-1,1}^{(2r)} = C_{-1,1}^{(2r+1)}  ,\ \ C_{-1,1}^{(1)}  = -1, \mbox{ and }
\]
\[
C_{-1,1}^{(2r+1)} = \left\{ \begin{array}{ll}
\frac{r+1}{2} & r \mbox{ is odd } \\
\frac{r-2}{2} & r \mbox{ is even}
\end{array}\right.
\]
\[
C_{-2}^{(r)} = \left\{\begin{array}{ll}
(-1,-1,1) & r\equiv 0 \mbox{ mod } 4\\
(1,1,-1) & r\equiv 1\mbox{ mod }4 \\
(0,-2,1) & r\equiv 2 \mbox{ mod }4 \\
(0,2,-1) & r\equiv 3 \mbox{ mod }4
\end{array}\right.
\]
\[
C_{-3}^{(r)} = C_{-1}^{(r+4)}
\]
\end{lem}

\subsection{Computation of degrees}
Now we are in position to prove the propositions \ref{degrees1} and \ref{degrees2}. The proposition \ref{degrees1} reduces to an easy check (using formulae for $C_{ij}'s$ from previous section) that degrees stated satisfy the following recurrence relations
\[
d_1^{(r+1)} = \left\{\begin{array}{ll}
d_1^{(r)} & r \mbox{ is odd}\\
-d_1^{(r)} + [C_{-1,1}^{(r)}]_- (1,-2,-1) + [C_{-2,1}^{(r)}]_-(3,0,0) + [C_{-3,1}^{(r)}]_-(1,2,1) & r\mbox{ is even}
\end{array}\right.
\]
\[
d_2^{(r+1)} = \left\{\begin{array}{ll}
d_2^{(r)} & r\mbox{ is odd }\\
-d_2^{(r)} + 3[C_{-1,2}^{(r)}]_-(1,-2,-1) + [C_{-2,1}^{(r)}]_-(3,0,0) + 3[C_{-3,2}^{(r)}]_-(1,2,1) & r\mbox{ is even}
\end{array}\right.
\]
\[
d_3^{(r+1)}= \left\{\begin{array}{ll}
d_3^{(r)} & r \mbox{ is even}\\
-d_3^{(r)} + [C_{-1,3}^{(r)}]_-(1,-2,-1)+[C_{-2,3}^{(r)}]_-(3,00)+[C_{-3,3}^{(r)}]_-(1,2,1) & r\mbox{ is odd}
\end{array}\right.
\]
together with initial values $d_1^{(0)} = (1,-1,-1)$, $d_2^{(0)} = (2,2,1)$ and $d_3^{(0)} = (3,3,2)$.\\

Similarly proposition \ref{degrees2} reduces to checking that degrees stated satisfy following recurrence relations:

\[
g_1^{(r+1)} = \left\{\begin{array}{ll}
g_1^{(r)} & r \mbox{ is odd}\\
-g_1^{(r)} + [C_{-1,1}^{(r)}]_- (1,-3,-2) + 3[C_{-2,1}^{(r)}]_-(1,0,0) + [C_{-3,1}^{(r)}]_-(1,3,2) & r\mbox{ is even}
\end{array}\right.
\]
\[
g_2^{(r+1)} = \left\{\begin{array}{ll}
g_2^{(r)} & r\mbox{ is odd }\\
-g_2^{(r)} + [C_{-1,2}^{(r)}]_-(1,-3,-2) + [C_{-2,1}^{(r)}]_-(1,0,0) + [C_{-3,2}^{(r)}]_-(1,3,2) & r\mbox{ is even}
\end{array}\right.
\]
\[
g_3^{(r+1)}= \left\{\begin{array}{ll}
g_3^{(r)} & r \mbox{ is even}\\
-g_3^{(r)} + [C_{-1,3}^{(r)}]_-(1,-3,-2)+3[C_{-2,3}^{(r)}]_-(1,0,0)+[C_{-3,3}^{(r)}]_-(1,3,2) & r\mbox{ is odd}
\end{array}\right.
\]
together with initial condition $g_1^{(0)} = (1,-3,-1)$, $g_2^{(0)} = (1,2,1)$ and $g_3^{(0)}=(2,3,2)$.

\section{Some determinant identities for $G_2$}\label{lemmas}
The aim of this part is to prove certain identities among the weight vectors of $V_{\omega_1}$ and $V_{\omega_2}$. 

\subsection{Determinant identities}
The proof of Lemmas \ref{seed1} and \ref{seed2}  follow from the following identities:
\begin{itemize}
\item[(D1)]
\[
 \Delta^{\omega_2}(\Delta^{s_2s_1\omega_1})^3 + (\Delta^{\omega_1})^3\Delta^{s_2s_1s_2\omega_2} = \Delta^{s_2\omega_2}\left( \Delta^{s_1\omega_1}(\Delta^{s_2s_1\omega_1})^2 - \Delta^{\omega_1}\Delta^{s_2s_1\omega_1}X_0 - (\Delta^{\omega_1})^2\Delta^{s_2s_1s_2s_1\omega_1}\right)
\]
\item[(D2)]
\[
 \begin{array}{ll}
(\Delta^{s_2s_1s_2s_1\omega_1})^3\Delta^{\omega_2} + \Delta^{w_0\omega_2}\left( \Delta^{s_1\omega_1}(\Delta^{s_2s_1\omega_1})^2 - \Delta^{\omega_1}\Delta^{s_2s_1\omega_1}X_0 - (\Delta^{\omega_1})^2\Delta^{s_2s_1s_2s_1\omega_1}\right)  = & \\
\Delta^{s_2s_1s_2\omega_2}\left(X_0(\Delta^{\omega_1}\Delta^{w_0\omega_1}+\Delta^{s_1\omega_1}\Delta^{s_2s_1s_2s_1\omega_1}) - \Delta^{\omega_1}\Delta^{s_1s_2s_1\omega_1}\Delta^{s_2s_1s_2s_1\omega_1} - \Delta^{s_1\omega_1}\Delta^{s_2s_1\omega_1}\Delta^{w_0\omega_1}\right) & 
\end{array}
\]
\item[(D3)]
\[
 \Delta^{s_1s_2s_1\omega_1}\Delta^{\omega_2} + \Delta^{\omega_1}\Delta^{s_1s_2\omega_2} = \Delta^{s_1\omega_1}F(1,1)
\]
\item[(D4)]
\[
 \Delta^{\omega_1}\Delta^{s_1s_2s_1s_2\omega_2} + \Delta^{w_0\omega_1}F(1,1) = \Delta^{s_1s_2s_1\omega_1}F^1(0,0)
\]
\end{itemize}

In order to prove these identities, one uses the fact that if $F\in V_{\lambda}(\mu)$ and $\mu<\lambda$, then $F=0$ if and only if $e_iF=0$ for $i=1,2$ (and using the action of operators $e_i's$ from Figures \ref{first} and \ref{second}).\\

For instance, consider $\xi = \Delta^{\omega_1}F(2,1)-\Delta^{s_1\omega_1}\Delta^{s_2\omega_2} + \Delta^{s_2s_1\omega_1}\Delta^{\omega_2}$. It is easy to compute and verify that $e_1\xi = e_2\xi = 0$ and hence $\xi=0$. More relations of this kind can be obtained by applying $f_1, f_2$ to $\xi$ (there are in fact 35 relations of this type, out of which 27 (=$dim(V_{2\omega_1})$) can be obtained from $\xi$). The relations (D1)-(D4) can be checked by direct calculation (for instance, (D3) is in fact obtained by applying $(f_1)^2$ to the equality $\xi=0$).

\subsection{Pl\"{u}cker relations}\label{plucker}
We will also need the following ``Pl\"{u}cker relations" which hold among weight vectors of $V_{\omega_2}$:

\begin{itemize}
\item[(P1)]
\[
F(2,1)^2 + \Delta^{\omega_2}F(1,0) = \Delta^{s_2\omega_2}F(1,1)
\]
\item[(P2)]
\[
F(1,0)\Delta^{s_1s_2s_1s_2\omega_2} = F^1(0,0)F(-2,-1)+F(1,1)\Delta^{w_0\omega_2}
\]
\item[(P3)]
\[
\Delta^{\omega_2}\Delta^{w_0\omega_2} = F(1,0)F(-1,0) - F^1(0,0)F^2(0,0)
\]
\item[(P4)]
\[
F(-2,-1)F(2,1) + F(1,0)F(-1,0) = F(1,1)F(-1,-1)
\]
\item[(P5)]
\[
F(1,1)^2 = F(2,1)\Delta^{s_1s_2\omega_2} - \Delta^{\omega_2}F(-1,0)
\]
\item[(P6)]
\[
\Delta^{\omega_2}\Delta^{s_1s_2s_1s_2\omega_2} = F(1,1)F(-1,0) - F^1(0,0)\Delta^{s_1s_2\omega_2}
\]
\end{itemize}

The proof of these relations is exactly same as the one in previous section. For example, the relation (P1) is equivalent to: $e_i^{\dagger}(F(2,1)^2 + \Delta^{\omega_1}F(1,0) - \Delta^{s_2\omega_2}F(1,1)) = 0$ for $i = 1,2$, which can be verified directly using the action of operators $e_i^{\dagger}$ from Figure \ref{second}. Other relations are obtained by applying lowering operators to both sides of (P1). For instance, (P5) can be obtained by applying $(1/4)(f_1^{\dagger})^2$ to (P1).

\section{Comparison with initial seed of \cite{gls}}\label{comparison}
In \cite{gls}, authors give a uniform combinatorial description of initial seeds for partial flag varieties. In this section we describe their construction for the case at hand. The algorithm in \cite{gls} first constructs the seed for $\C[N_i]$ (where $N_i$ is unipotent radical of $P_i$ for $i=1,2$) and then lift it to a seed of $\C[P_j^-\backslash G]$ ($j = \{1,2\}\setminus \{i\}$). In the next two paragraphs we will give this lifting. Reader can refer to \cite{gls} for details.

\subsection{Case $i=1$} In this case the initial seed for $\C[P_2^-\backslash G]$, denoted by $\Sigma_1^{(GLS)}$ has exchange matrix given by
\[
\Gamma_1^{(GLS)}:
\xymatrix{
{\bf -1} \ar[r] & 1 \ar[r] \ar[d] & 3 \ar[d] & {\bf -3} \ar[l] \\
 & {\bf -2} \ar[r] & 2 \ar[lu] & 
}
\]
with corresponding functions given by 
\[
X_{-1}^{(GLS)} = \Delta^{w_0\omega_1},\ X_{-3}^{(GLS)} = \Delta^{\omega_1}
\]
\[
X_{-2}^{(GLS)} = X_0(\Delta^{\omega_1}\Delta^{w_0\omega_1}+\Delta^{s_1\omega_1}\Delta^{s_2s_1s_2s_1\omega_1}) - \Delta^{\omega_1}\Delta^{s_1s_2s_1\omega_1}\Delta^{s_2s_1s_2s_1\omega_1} - \Delta^{s_1\omega_1}\Delta^{s_2s_1\omega_1}\Delta^{w_0\omega_1}
\]
\[
X_1^{(GLS)} = -\Delta^{\omega_1}\Delta^{s_1s_2s_1\omega_1} + \Delta^{s_1\omega_1}X_0
\]
\[
X_2^{(GLS)} = \Delta^{s_2s_1\omega_1}(\Delta^{s_1\omega_1})^2 - 2\Delta^{\omega_1}\Delta^{s_1\omega_1}X_0 + (\Delta^{\omega_1})^2\Delta^{s_1s_2s_1\omega_1},\ X_3^{(GLS)} = \Delta^{s_1\omega_1}
\]

\noindent {\bf Claim:} The sequence of mutations $\mu_1\mu_3$ applied to $\Sigma_1^{(GLS)}$ produces same cluster as $\underline{\Sigma_1}^0$ with exchange matrix $-\underline{\widetilde{B_1}}^0$.\\

Since $\underline{\Sigma_1}^0$ is obtained from initial seed $\underline{\Sigma_1}$ by applying mutations $\mu_2\mu_3$, in order to prove the claim we need to compute $\mu_3\mu_2\mu_1\mu_3(\Sigma_1^{(GLS)})$. We do this computation in following steps:

\begin{itemize}
\item[a)] $\mu_3(\Gamma_1^{(GLS)}$ is:
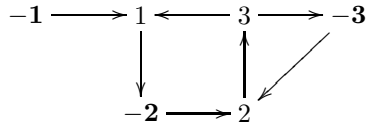
\[
\xymatrix{
{\bf -1} \ar[r] & 1 \ar[d] & 3 \ar[l] \ar[r] & {\bf -3} \ar[ld] \\
 & {\bf -2} \ar[r] & 2 \ar[u] &
}
\] 
And the cluster variable at third vertex mutates as:
\[
X_3' = \frac{\Delta^{\omega_1}(\Delta^{s_1\omega_1}X_0 - \Delta^{\omega_1}\Delta^{s_1s_2s_1\omega_1}) + (\Delta^{s_1\omega_1})^2\Delta^{s_2s_1\omega_1} - 2\Delta^{s_1\omega_1}\Delta^{\omega_1}X_0+(\Delta^{\omega_1})^2\Delta^{s_1s_2s_1\omega_1}}{\Delta^{s_1\omega_1}}
\]
\[
X_3' = \Delta^{s_1\omega_1}\Delta^{s_2s_1\omega_1} - \Delta^{\omega_1}X_0
\]

\item[b)] $\mu_1\mu_3(\Gamma_1^{(GLS)})$ is:
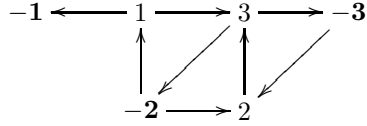
\[
\xymatrix{
{\bf -1} & 1 \ar[r] \ar[l] & 3 \ar[ld] \ar[r] & {\bf -3} \ar[ld] \\
 & {\bf -2} \ar[u] \ar[r] & 2 \ar[u] &
}
\]
And the cluster variable at vertex $1$ mutates as:
\[
X_1' = \frac{X_{-1}X_3' + X_{-2} }{X_1 }
\]
This exchange relation gives $X_1' = \Delta^{s_2s_1s_2s_1\omega_1}$.

\item[c)] $\mu_2\mu_1\mu_3(\Sigma_1^{(GLS)})$ is:
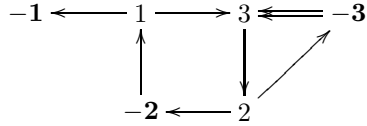
\[
\xymatrix{
{\bf -1} & 1 \ar[r] \ar[l] & 3 \ar[d] & {\bf -3} \ar[l]<1pt> \ar[l]<-1pt> \\
& {\bf -2} \ar[u] & 2 \ar[l] \ar[ru]
}
\]
The new cluster variable at vertex $2$ is obtained as:
\[
X_2' = \frac{(X_3')^3 + X_{-2}X_{-1}}{X_2}
\]
This can be easily computed to be $X_2' = \Delta^{s_1\omega_1}(\Delta^{s_2s_1\omega_1})^2 - \Delta^{\omega_1}\Delta^{s_2s_1\omega_1}X_0 - (\Delta^{\omega_1})^2\Delta^{s_2s_1s_2s_1\omega_1}$.

\item[d)] Finally $\mu_3\mu_2\mu_1\mu_3(\Gamma_1^{(GLS)})$ is:
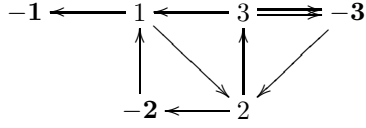
\[
\xymatrix{
{\bf -1} & 1 \ar[l] \ar[rd] & 3 \ar[l] \ar[r]<1pt> \ar[r]<-1pt> & {\bf -3} \ar[ld] \\
& {\bf -2} \ar[u] & 2 \ar[l] \ar[u] &
}
\]
which is same as our initial exchange matrix except for all arrows reversed. And the last cluster variable here (at vertex $3$) is:
\[
X_3'' = \frac{(\Delta^{\omega_1})^2\Delta^{s_2s_1s_2s_1\omega_1}+\Delta^{s_1\omega_1}(\Delta^{s_2s_1\omega_1})^2 - \Delta^{\omega_1}\Delta^{s_2s_1\omega_1}X_0 - (\Delta^{\omega_1})^2\Delta^{s_2s_1s_2s_1\omega_1} }{\Delta^{s_1\omega_1}\Delta^{s_2s_1\omega_1} - \Delta^{\omega_1}X_0}
\]
Therefore $X_3'' = \Delta^{s_2s_1\omega_1}$.
\end{itemize}

Thus we have shown that our initial seed $\underline{\Sigma_1}$ can be obtained from the one given in \cite{gls} by a sequence of mutations (up to a sign). Since the cluster algebras associated to matrix $\widetilde{B}$ and $-\widetilde{B}$ are naturally isomorphic, both initial seeds give rise to same cluster algebra.

\subsection{Case $i=2$} We start with reduced word ${\bf i_1} = (1,2,1,2,1,2)$ of $w_0$. The initial seed for $\C[P_1^-\backslash G]$, denoted by $\Sigma_2^{(GLS)}$ has exchange matrix given by:
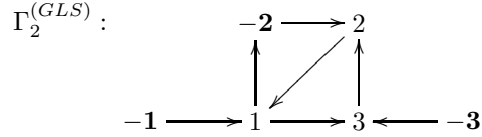
\[
\Gamma_2^{(GLS)}:
\xymatrix{
 & {\bf -2} \ar[r] & 2 \ar[ld] & \\
{\bf -1} \ar[r]& 1 \ar[u] \ar[r] & 3 \ar[u] & {\bf -3} \ar[l]
}
\]
with corresponding functions given by
\[
Y_{-1}^{(GLS)} = \Delta^{w_0\omega_2},\ Y_{-2}^{(GLS)} = F^1(0,0),\ Y_{-3}^{(GLS)} = \Delta^{\omega_2}
\]
\[
Y_1^{(GLS)} = F(2,1)F(1,0) - \Delta^{s_2\omega_2}F^1(0,0),\ Y_2^{(GLS)} = F(2,1),\ Y_3^{(GLS)} = \Delta^{s_2\omega_2}
\]

\noindent {\bf Claim:} The sequence of mutations $\mu_1\mu_3$ applied to $\Sigma_2^{(GLS)}$ yields cluster of $\underline{\Sigma_2}^0$ with exchange matrix multiplied by $-1$.\\

The proof of similar claim from last paragraph carries verbatim over to proof of this claim. One only needs to make use of the ``Pl\"{u}cker Relations" ((P1)-(P6)) from last section in order to carry out computations involving exchange relations.

\section{Summary of results}
In this section we give brief summary of the results proved in this note.

\subsection{Case $i=1$}
\begin{itemize}
\item[-] The initial seed $\underline{\Sigma_1}$ can be obtained from initial seed $\Sigma_1$ of $\A$ by applying mutations $\mu_2\mu_4$ and freezing vertex $2$ (see section \ref{explanation}).
\item[-] The cluster algebra associated to initial seed $\underline{\Sigma_1}$ is isomorphic to the algebra $\A_1$ localized at multiplicative set $\{(\Delta^{\omega_1})^n\}_{n\in \Z_{\geq 0}}$. Explicitly we obtain all the weight vectors of Figure \ref{first} by applying mutations to $\underline{\Sigma_1}$, as follows (see section \ref{proof1}):
\[
\Delta^{\omega_1} = X_{-3},\ \ \Delta^{w_0\omega_1} = X_{-1},\ \ \Delta^{s_2s_1\omega_1} = X_3,\ \ \Delta^{s_2s_1s_2s_1\omega_1} = X_1
\]
\[
\Delta^{s_1\omega_1} = X_1^{(-2)} = \mbox{ function at vertex 1 of } \mu_1\mu_3\mu_2\mu_3(\underline{\Sigma_1})
\]
\[
\Delta^{s_1s_2s_1\omega_1} = X_1^{(-6)} = \mbox{ function at vertex 1 of } \mu_1\mu_3(\mu_1\mu_2\mu_3)^2\mu_2\mu_3(\underline{\Sigma_1})
\]
And finally we have $\Delta^{\omega_1}X_0 = \Delta^{s_1\omega_1}\Delta^{s_2s_1\omega_1}-X_3^{(0)}$.

\item[-] The initial seed $\underline{\Sigma_1}$ can also be obtained by applying sequence of mutations $\mu_3\mu_2\mu_1\mu_3$ to initial seed $\Sigma_1^{(GLS)}$ as given in \cite{gls} (see section \ref{comparison}).
\end{itemize}

\subsection{Case $i=2$} 
\begin{itemize}
\item[-] The initial seed $\underline{\Sigma_2}$ can be obtained from the initial seed $\Sigma_2$ of $\A$ by applying mutations $\mu_2\mu_4$ and freezing vertex $2$ (section \ref{explanation}).
\item[-] The cluster algebra associated to initial seed $\underline{\Sigma_2}$ is isomorphic to $\A_2$ localized at the multiplicative set $\{(\Delta^{\omega_2})^n\}_{n\in \Z_{\geq 0}}$. Explicitly we obtain weight vectors of Figure \ref{second} as (see section \ref{proof2}):
\[
\Delta^{\omega_2} = Y_{-3},\ \ \Delta^{w_0\omega_2} = Y_{-1},\ \ F^1(0,0) = Y_{-2}
\]
\[
\Delta^{s_1s_2s_1s_2\omega_2} = Y_1,\ \ F(1,1) = Y_2,\ \ \Delta^{s_1s_2\omega_2} = Y_3
\]
\[
\Delta^{s_2\omega_2} = Y_1^{(-2)}
\]
\[
F(2,1) = Y_2^{(0)}
\]
\[
\Delta^{s_2s_1s_2\omega_2} = Y_1^{(-6)}
\]
\[
F(-1,-1) = Y_2^{(-4)}
\]
\[
F(-2,-1) = Y_2^{(-2)}
\]
\[
F(1,0) \mbox{ and } F(-1,0) \mbox{ are functions at vertex 2 of } \mu_2\mu_3\mu_1(\underline{\Sigma_2}^0) \mbox{ and } \mu_2\mu_3\mu_2(\underline{\Sigma_2}^0) \mbox{ respectively}
\]
Finally we have following Pl\"{u}cker relation:
\[
F^2(0,0)\Delta^{\omega_2} = \Delta^{s_2\omega_2}\Delta^{s_1s_2\omega_2} - F(1,1)F(2,1) - \Delta^{\omega_2}F^1(0,0)
\]

\item[-] The initial seed $\underline{\Sigma_2}$ can also be obtained by applying $\mu_3\mu_2\mu_1\mu_3$ to the initial seed $\Sigma_2^{(GLS)}$ (see section \ref{comparison}).
\end{itemize}

\begin{center}
{\sc Acknowledgments}
\end{center}
I would like to thank Prof. Andrei Zelevinsky for introducing me to the theory of cluster algebras and for several helpful and stimulating discussions. I am also grateful to Prof. Bernard Leclerc for kindly explaining to me the construction of initial seeds, as in \cite{gls}. \\


\begin{thebibliography}{10}
\bibitem{schubert} A. Berenstein, A. Zelevinsky {\em Total positivity in Schubert varieties} Comment. Math. Helv. {\bf 72} (1997) No. 1, 128-166.
\bibitem{CA3} A. Berenstein, S. Fomin, A. Zelevinsky {\em Cluster algebras III: Upper bounds and double Bruhat cells} Duke Math. J. {\bf 126} (2005) No. 1, 1-52.
\bibitem{giovanni} Giovanni Cerulli Irelli {\em Structural theory of rank three cluster algebras of affine type} Ph.D. Thesis at Universit\`{a} degli studi di Padova, Italy.
\bibitem{demonet} L. Demonet  {\em Cluster algebras and preprojective algebras: the non simply-laced case} arxiv:0711.4098v1
\bibitem{double} S. Fomin, A. Zelevinsky {\em Double Bruhat cells and Total positivity} J. Amer. Math. Soc. {\bf 12} (1999) No. 2, 335-380.
\bibitem{CA1} S. Fomin, A. Zelevinsky {\em Cluster algebras I: Foundations} J. Amer. Math. Soc. {\bf 15} (2002) No. 2, 497-529.
\bibitem{CA4} S. Fomin, A. Zelevinsky {\em Cluster algebras IV: Coefficients} Compositio Mathematica {\bf 143} (2007) 112-164.
\bibitem{gls} C. Geiss, B. Leclerc, J. Schr\"{o}er {\em Partial flag varieties and Preprojective algebras} arXiv:math/0609138
\end{thebibliography}
\end{document}